\documentclass{article}
\usepackage{amsmath}
\usepackage{amssymb}
\usepackage{amsthm}
\usepackage[bookmarks,colorlinks=false]{hyperref}
\usepackage{comment}

\newtheorem{theorem}{Theorem}
\newtheorem{definition}{Definition}
\newtheorem{lemma}{Lemma}
\newtheorem{remark}{Remark}

\newtheorem{proposition}{Proposition}

\newcounter{mnotecount}[section]

\newcommand{\mnotex}[1]
{\protect{\stepcounter{mnotecount}}$^{\mbox{\footnotesize $\bullet$\themnotecount}}$
\marginpar{
\raggedright\tiny\em
$\!\!\!\!\!\!\,\bullet$\themnotecount: #1} }

\begin{document}

\title{\textbf{The $\mathcal{C}$-connection and the 4-dimensional
Einstein spaces}}

\author{ A. Garc\'\i a-Parrado \\
	Departamento de Matem\'aticas, Universidad de C\'ordoba,\\
	Campus de Rabanales, 14071, C\'ordoba, Spain\\
	e-mail: agparrado@uco.es  \\[2ex]
}
\maketitle

\begin{abstract}
We introduce a new family of operators
in 4-dimensional pseudo-Riemannian manifolds
with a non-vanishing Weyl scalar (non-degenerate spaces)
that keep the conformal covariance of
\emph{conformally covariant tensor concomitants}. A
particular case that arises naturally is the $\mathcal{C}$-connection
that is a Weyl connection that keeps \emph{conformal invariance}.
Using the $\mathcal{C}-$connection we give a new characterization of
non-degenerate spaces that are conformal to an Einstein space.

\end{abstract}


\section{Introduction}

Conformal geometry has become a crucial area of study in both physics and mathematics due to its role in
describing fundamental interactions. All fundamental interactions in physics, except for gravity,
are governed by \emph{conformally covariant} dynamical equations.
This naturally leads to the important mathematical question of
identifying and classifying geometric objects that are \emph{conformally covariant}. Despite the significant
attention this question has received, there is no unified theory of conformal covariance. Instead, research has
progressed through targeted approaches, with a variety of distinct perspectives and applications.

A comprehensive introduction to conformal geometry is provided in \cite{CURRYGOVER},
where the role of  \emph{tractor calculus} is emphasized.
Another interesting outcome is \cite{szekeres1968conformal} where \emph{conformal tensors}
are classified and their properties as metric concomitants are fully investigated.

Higher-order differential operators, such as the Paneitz operator \cite{paneitz1983quartic}, have played a
significant role in conformal geometry. The Paneitz operator, a fourth-order conformally invariant operator,
generalizes the Laplacian and is crucial in geometric analysis. Similarly,
\cite{graham1992conformally} provides a geometric construction of conformally invariant powers of the Laplacian,
particularly using the Fefferman-Graham ambient metric. Their work highlights the existence of conformally
invariant operators in various dimensions.

Other relevant contributions are \cite{fialkow1944conformal}, on the conformal differential
geometry of subspaces and the book \cite{juhl2009families}, exploring the connections between conformal
geometry, higher-order differential operators, and physical theories such as holography.

In general relativity, the introduction of the concept of the \emph{conformal
boundary} by Penrose in the 1960s \cite{Pen63} marked the
birth of a new field with significant impact on theoretical physics and pure
mathematics. In physics, this concept
is vital for studying isolated systems and their radiative properties, while in
mathematics, it helps define the
causal boundary of a Lorentzian manifold.

Another major advancement in general relativity
came with H. Friedrich's work on the \emph{conformal equations} \cite{Fri81a,
Fri81b}.
Friedrich's equations have been extensively studied (see \cite{Fri02, JUANBOOK}) and have led to significant
applications \cite{Fri86a, Fri86b, Fri91}.

Still within the framework of general relativity, an important question,
also addressed in the Friedrich formalism is to determine the
class of space-times that are conformal to a vacuum solution of the
Einstein field equations. A number of approaches to answer this question
can be found in the literature
\cite{KoNewTod85,Listing2001,Edgar2004,GOVER2006450,szekeres1963spaces}.

In this work we perform an analysis of conformally
covariant \emph{metric concomitants} defined in four-dimensional
pseudo-Riemannian manifolds (see section \ref{sec:conformal-pseudo-differential}).
If the pseudo-Riemannian manifold is \emph{non-degenerate} (the
Weyl scalar $C_{abcd}C^{abcd}$ is no-where vanishing) then
we can define the notion of \emph{conformal pseudo-differential} (see
eqns. \eqref{eq:define-conformal-diff-scalars},
\eqref{eq:dspq-diff}). Using
this tool we can construct new conformally covariant metric concomitants
out of a given one as explained in
subsections \ref{subsec:conf-diff-scalars} and
\ref{subsec:conf-diff-tensors}. In the particular case of a
\emph{conformally invariant tensor concomitant} we can introduce the
$\mathcal{C}-$connection (see subsection \ref{subsec:C-connection}).
This is a Weyl connection that keeps the conformal invariance
as stated by Theorem \ref{theo:c-connection-conf-invariant}. In section
\ref{sec:conformal-einstein} we present an application consisting
in a new characterization of non-degenerate pseudo-Riemannian manifolds
conformal to Einstein spaces (see Theorem \ref{theo:ideal-einstein}).
This characterization is \emph{algorithmic} in the
sense that once a pseudo-Riemannian metric is given it can be tested whether
it is conformal to an Einstein space by a procedure with a finite number of
steps.

All the computations of this paper have been carried out with \cite{XACT}
(see also \cite{XPERM}).

\section{Geometric preliminaries}
\label{sec:preliminaries}
Let $\mathcal{M}$ be a 4-dimensional smooth manifold (unless otherwise stated, all structures are assumed to be smooth) and suppose that we define on $\mathcal{M}$ a
pseudo-Riemannian metric $\tilde{g}_{ab}$ (referred to as the physical metric) and
a second pseudo-Riemannian metric, $g_{ab}$ (called the unphysical metric) of the same signature
which is {\em conformally related} to the first in the following fashion
\begin{equation}
g_{ab} = \Theta^2 \tilde{g}_{ab}.
\label{eq:unphysicaltophysical-downstairs}
\end{equation}
The conformal factor $\Theta$ is assumed to be a positive smooth function which does not vanish
on $\mathcal{M}$.
The previous relation can be used to define two pseudo-Riemannian spaces:
$(\mathcal{M}, g_{ab})$ and $(\mathcal{M}, \tilde{g}_{ab})$.
The conformal relation \eqref{eq:unphysicaltophysical-downstairs}
can be used to set up an equivalence relation between
conformally related pseudo-Riemannian metrics. The resulting equivalence classes are called
\emph{conformal structures} on $\mathcal{M}$.
In this way, if \eqref{eq:unphysicaltophysical-downstairs}
holds then $g_{ab}$ and $\tilde{g}_{ab}$ belong to the same conformal structure.

The bundle of $r$-contravariant and $s$-covariant tensors
$T^r_s(\mathcal{M})$ is introduced in the standard way
(the union of all these bundles will be denoted by $T(\mathcal{M})$).
Tensor fields (or just tensors) are sections of these bundles and
the corresponding set of sections is ${\mathcal T}^r_s(\mathcal{M})$,
${\mathcal T}(\mathcal{M})$.
We shall use small Latin letters to denote abstract indices
of tensors in $\mathcal{T}(\mathcal{M})$.
Indices are always raised and lowered with respect to the unphysical metric $g_{ab}$
with the exception of $\tilde{g}^{ab}$, where we follow
the traditional convention that it represents the inverse
of $\tilde{g}_{ab}$
\begin{equation}
\tilde{g}^{ac}\tilde{g}_{cb}=\delta_b{}^a.
\end{equation}
From \eqref{eq:unphysicaltophysical-downstairs}, the explicit relation between
the physical and the unphysical contravariant metric tensors becomes
\begin{equation}
g^{ab} = \frac{\tilde{g}^{ab}}{\Theta^2}.
\label{eq:unphysicaltophysical-upstairs}
\end{equation}

Also each metric tensor has its own Levi-Civita connection denoted respectively
by $\nabla_a$, $\tilde{\nabla}_a$ which are used to define the connection coefficients and the curvature
tensors in the standard fashion. Our conventions for the (unphysical) Riemann, Ricci and Weyl tensors
are
\begin{equation}
\nabla_{a}\nabla_{b}\omega_{c} -  \nabla_{b}\nabla_{a}\omega_{c} =  R_{abc}{}^{d} \omega_{d}\;,
\label{eq:define-riemann}
\end{equation}
\begin{equation}
 R_{ac}\equiv R_{abc}{}^b\;,
\label{eq:define-ricci}
\end{equation}
\begin{equation}
C_{abc}{}^{d}\equiv R_{abc}{}^{d} - 2 L^{d}{}_{[b}g_{a]c} - 2\delta^{d}{}_{[b}L_{a]c}\;,
\label{eq:define-weyl}
\end{equation}
where the unphysical Schouten tensor is defined by
\begin{equation}
 L_{ab} \equiv \tfrac{1}{2} (R_{ab} -  \tfrac{1}{6} R g_{ab}).
\label{eq:define-schouten}
\end{equation}
Tensors defined in terms of the physical metric $\tilde{g}_{ab}$ will be denoted with a tilde over
the symbol employed for an unphysical spacetime tensor.
The difference between the Levi-Civita connection
coefficients of the physical and the unphysical metrics
is the \emph{transition tensor}
$\Gamma[\nabla ,\tilde{\nabla}]^{e}{}_{ac}$ given by
\begin{equation}
\Gamma[\nabla ,\tilde{\nabla}]^{e}{}_{ac} =
2 \delta_{(a}{}^{e}\Upsilon_{c)} -  g_{ca} \Upsilon^{e},
\label{eq:define-connection-tensor}
\end{equation}
where we defined
\begin{equation}
\Upsilon_a\equiv\frac{\nabla_a\Theta}{\Theta}.
\label{eq:define-upsilon}
\end{equation}

Any (tensorial) object $Q$ defined exclusively in terms of the unphysical metric
$g_{ab}$ and its Levi-Civita connection $\nabla$ will be
called an \emph{unphysical metric concomitant}.
Elementary examples
of metric concomitants are the Riemann curvature of the Levi-Civita connection,
its non-trivial metric contractions (Ricci tensor and Ricci scalar) and the
Weyl and Schouten tensors.

Of course, any unphysical metric concomitant will always admit a physical
counterpart (\emph{physical metric concomitant}),
denoted with a tilde on the top, $\tilde{Q}$. We will just speak of a
metric concomitant if the context makes it clear what the involved metric is.
Sometimes we may stress the dependencies on the metric and the Levi-Civita
connection in the notation (such as in Definition \ref{def:conf-metric-concomitant}).

\section{Conformal pseudo-differential
of conformally covariant metric concomitants}
\label{sec:conformal-pseudo-differential}

\begin{definition}[Conformally covariant metric concomitant]
We will say that a metric concomitant $Q(g,\nabla)$ is
conformally covariant of weight $s\in\mathbb{Z}$ if
$Q(g,\nabla)$ and its physical
counterpart $\tilde Q(\tilde g,\tilde \nabla)$ are related to each other
in the following way
\begin{equation}
 \tilde Q(\tilde g,\tilde \nabla) = \Theta^s Q(g,\nabla),
 \quad s\in\mathbb{Z}.
\end{equation}
\label{def:conf-metric-concomitant}
\end{definition}
Classical examples of conformally covariant metric concomitants are the metric tensor
itself and its inverse (see eqs. \eqref{eq:unphysicaltophysical-downstairs},
\eqref{eq:unphysicaltophysical-upstairs}), the Weyl tensor
$\tilde{C}_{abc}{}^{d}=C_{abc}{}^{d}$ and the Bach tensor
\begin{equation}
 B_{bc}\equiv \nabla^a\nabla^d C_{abdc}-\frac{1}{2}R^{ad}C_{abcd},
 \label{eq:bach-tensor}
\end{equation}
which is conformally invariant $\tilde{B}_{ab}(\tilde g,\tilde\nabla)
=B_{ab}(g, \nabla)$.
See \cite{szekeres1968conformal} for an in-depth study of
conformally covariant metric concomitants in pseudo-Riemannian
manifolds of arbitrary dimension.

Define the following metric concomitant
$C\cdot C\equiv C_{abcd}C^{abcd}$.
Unless otherwise stated we will assume that
$(\mathcal{M},g_{ab})$ fulfills the condition $(C \cdot C)\neq 0$
at every point.
Let us define the following unphysical metric concomitant:
\begin{equation}
\Lambda_d = \frac{8 C_{d}{}^{amp}\nabla_{[m} L_{p]a}}{(C \cdot C)}.
\label{eq:define-Lambda}
\end{equation}

The following lemma is key for our computations.
\begin{lemma}
The unphysical $\Lambda_d$ just defined and its physical
counterpart $\tilde{\Lambda}_d$ fulfill the relation
\begin{equation}
\Upsilon_a = \Lambda_a - \widetilde{\Lambda}_a.
\end{equation}
\label{lemm:lambda-upsilon}
\end{lemma}
\proof

The relation between the
unphysical and physical Schouten tensors can be always written in the form (see e.g. \cite{JUANBOOK})
\begin{equation}
L_{ab} = \tilde{L}_{ab} -  \Upsilon_{a} \Upsilon_{b} + \frac{1}{2} g_{ab} \Upsilon_{c} \Upsilon^{c} -  \nabla_{b}\Upsilon_{a}.
\label{eq:nabla-upsilon}
\end{equation}
Computing the unphysical covariant derivative of the previous equation and using \eqref{eq:define-riemann} with $\omega_d$ replaced by
$\Upsilon_d$ we get the following integrability condition
\begin{equation}
2\tilde{L}_{b[p} g_{m]a} \Upsilon^{b}
+2\tilde{L}_{a[m}\Upsilon_{p]} +
\Upsilon^{b} C_{abmp} + 2\nabla_{[p}\tilde{L}_{m]a} +
2\nabla_{[m}L_{p]a} = 0.
\label{eq:integrability-schouten}
\end{equation}
Expressing $\nabla$ in terms of $\tilde\nabla$ by means of \eqref{eq:define-upsilon}
we obtain the following equivalent form of \eqref{eq:integrability-schouten}
\begin{equation}
\tilde{\nabla}_{[b}\tilde{L}_{e]a} = \tfrac{1}{2} \Upsilon^{c} C_{acbe} + \nabla_{[b}L_{e]a}.
\label{eq:diff-tildeL-L}
\end{equation}
On the other hand the relation between the
unphysical Weyl tensor $C_{abcd}$ and its
physical counterpart, denoted by $\tilde{W}_{abcd}$, is
\begin{equation}
C_{abcd} = \Theta^2 \tilde{W}_{abcd}.
\label{eq:weylunphys-cov-to-weylphys}
\end{equation}
Now, take eq. \eqref{eq:diff-tildeL-L}, multiply both sides by
$\tilde{g}^{ap}\tilde{g}^{bq}\tilde{g}^{er}\tilde{W}_{dpqr}$, use \eqref{eq:weylunphys-cov-to-weylphys}
and the following identity, valid only in dimension four,
(see e.g. \cite{EDGARHOGLUNDDDI}),
\begin{equation}
C_{abmp}C^{dbmp}=\tfrac{1}{4}\delta_a{}^dC_{cbmp}C^{cbmp}.
\end{equation}
The sought result follows after elementary algebraic manipulations.
\qed

\medskip

\subsection{Conformal differentiation of scalar
conformally co\-va\-riant metric concomitants}
\label{subsec:conf-diff-scalars}
\noindent
Let us consider the particular case of a conformally covariant
metric concomitant defined by a scalar function $w\in C^\infty(\mathcal{M})$.
In this case and following the notation of Definition \ref{def:conf-metric-concomitant}
we have
\begin{equation}
\widetilde{w} = \Theta^s w.
\end{equation}
Taking the $\tilde\nabla$ covariant derivative of this
equation and replacing it by $\nabla$ where appropriate
we get
\begin{equation}
 \widetilde{\nabla}_a \widetilde{w} = s \Theta^{s-1} \nabla_a \Theta w + \Theta^s \nabla_a w.
 \label{eq:cd-not-conf}
\end{equation}
This equation shows that the Levi-Civita connection is not
\emph{conformally covariant}. However, we can use it as the starting
point to construct a \emph{conformally covariant operator}.
Write \eqref{eq:cd-not-conf} in terms of $\Upsilon_a$
and use Lemma \ref{lemm:lambda-upsilon}
\[
\widetilde{\nabla}_a\widetilde{w} = \Theta^s\nabla_a w + s\Theta^s\Upsilon_a w,
\]

\[
\Leftrightarrow \widetilde{\nabla}_a \widetilde{w} = s \Theta^{s} (\Lambda_a - \widetilde{\Lambda}_a) w + \Theta^s \nabla_a w,
\]

\[
\Leftrightarrow \widetilde{\nabla}_a \widetilde{w} + s \Theta^{s} \widetilde{\Lambda}_a w = \Theta^s (\nabla_a w + s \Lambda_a w),
\]

The previous equation can be written in terms of
the following operator $D^s_a:C^\infty(\mathcal{M})\rightarrow
C^\infty(\mathcal{M})$
\begin{equation}
\quad D^s_a w \equiv \nabla_a w + s \Lambda_a w.
\label{eq:define-conformal-diff-scalars}
\end{equation}
The previous considerations lead us to the proof
of the following result.

\begin{theorem}
The operator $D^s_a$ is conformally covariant when acting on
conformally covariant scalars.
\begin{equation}
\widetilde{D}^s_a \widetilde{w} = \Theta^s D^s_a w.
\end{equation}
\label{theo:conformal-diff-scalars}
 \end{theorem}
\qed

Another straightforward property of $D^s_a$
is contained in the following proposition.

\begin{proposition}[Generalized Leibnitz rule]
The operator
$D_a^{s}:C^\infty(\mathcal{M})\rightarrow C^\infty(\mathcal{M})$
is linear and fulfills the following property (generalized Leibnitz rule)
\begin{equation}
 D^{s}_a(w_1 w_2) = (D^{s_1}_a w_1) w_2+
 w_1 (D_a^{s_2} w_2),\quad s = s_1+s_2,\quad
 \forall w_1, w_2\in C^{\infty}(\mathcal{M}).
 \label{eq:leibnitz}
\end{equation}
\end{proposition}
\proof
The operator $D^s_a$ is obviously linear, so we only need
to show that it verifies \eqref{eq:leibnitz}. To that end we work out
its left hand side and obtain its right hand side.
\begin{eqnarray*}
 &&
 D^{s}_a(w_1 w_2) = \nabla_a (w_1 w_2) + s \Lambda_a (w_1 w_2) =\nonumber\\
 &&
w_2 \nabla_a w_1 + w_1 \nabla_a w_2+ s_1 \Lambda_a w_1 w_2
+s_2 \Lambda_a w_1 w_2=\nonumber\\
&&
w_1(\nabla_a w_2+s_2\Lambda_a  w_2)+
w_2(\nabla_a w_1+s_1\Lambda_a w_1)=
(D^{s_1}_a w_1) w_2+ w_1 (D_a^{s_2} w_2).
\end{eqnarray*}

\qed

\subsection{Conformal pseudo-differential of conformally cova\-riant tensors}
\label{subsec:conf-diff-tensors}
Consider now a metric concomitant that is a conformally covariant tensor of rank $p$
and conformal weight $s$. Then, according to definition \ref{def:conf-metric-concomitant} one has
\begin{equation}
 \tilde{K}_{a_1 a_{2} \ldots a_{p}} = \Theta^{s} K_{a_1 a_{2} \ldots a_{p}}.
\label{eq:conformal-covariant-K}
\end{equation}

Take the physical covariant derivative of eq. \eqref{eq:conformal-covariant-K}
\[
\tilde{\nabla}_{b} \tilde{K}_{a_1 a_{2} \ldots a_{p}} = s \Theta^{s-1} \nabla_{b}
\Theta K_{a_1 a_{2} \ldots a_{p}} + \Theta^{s}
\tilde{\nabla}_{b} K_{a_1 a_{2}
\ldots a_{p}}.
\]
Next, use the transition tensor \eqref{eq:define-connection-tensor}
on the right hand side to write $\tilde\nabla$
in terms of $\nabla$
\[
s \Theta^{s-1} \nabla_{b} \Theta K_{a_1 a_{2} \ldots a_{p}} + \Theta^{s} \left(
\nabla_{b} K_{a_1 a_{2} \ldots a_{p}} + \sum_{j=1}^p \Gamma[\nabla ,\tilde{\nabla}]^{c_{j}}{}_{b a_j}
K_{a_{1} a_{2} \ldots c_{j} \ldots a_{p}} \right)
\]

\begin{equation}
= \Theta^{s} \left[ \sum_{j=1}^p M^{c_{j}}{}_{b a_{j}}({\boldsymbol \Upsilon},{\boldsymbol g})
K_{a_{1} \ldots c_{j} \ldots a_{p}} + \nabla_{b} K_{a_{1} \ldots a_{p}} \right],
\label{eq:tildenablaK-transition}
\end{equation}
where in the last step we defined (we use boldface to mean that the covariant
form of a tensor is understood)
\begin{equation}
M^{c}{}_{b a}({\boldsymbol \Upsilon},{\boldsymbol g}, p, s)
\equiv \frac{s}{p} \Upsilon_{b} \delta_{a}{}^{c} + \Gamma[\nabla ,\tilde{\nabla}]^c{}_{b a}.
\end{equation}

Expanding out the transition tensor (see \eqref{eq:define-connection-tensor})
we get the equivalent form
\[
M^{c}{}_{b a}({\boldsymbol \Upsilon},{\boldsymbol g}, p, s) =
\frac{s}{p} \Upsilon_{b} \delta_{a}{}^{c} + \Upsilon_{a} \delta_{b}{}^{c}
+ \Upsilon_{b} \delta_{a}{}^{c} - g_{a b} g^{c r} \Upsilon_{r}
\]

\[
= \left( \frac{s + p}{p} \right) \Upsilon_{b} \delta_{a}{}^{c} + \Upsilon_{a} \delta_{b}{}^{c} - g_{a b} g^{c r} \Upsilon_{r}
\]
Note that $M^{c}{}_{b a}({\boldsymbol \Upsilon},{\boldsymbol g},p,s)$ is
linear in its first argument and conformally invariant in the second one.
Hence using the relation \( \Upsilon_{a} = \Lambda_{a} - \tilde{\Lambda}_{a}\)
proven by Lemma \ref{lemm:lambda-upsilon}
we get:
\begin{equation}
 M^{c}{}_{b a}({\boldsymbol \Upsilon},{\boldsymbol g},p,s)=
 M^{c}{}_{b a}({\boldsymbol\Lambda} - \tilde{\boldsymbol\Lambda},{\boldsymbol g},p,s)=
 M^{c}{}_{b a}({\boldsymbol\Lambda},{\boldsymbol g},p,s)
 - M^{c}{}_{b a}(\tilde{\boldsymbol\Lambda},\tilde{\boldsymbol g},p,s).
\end{equation}

Putting this back in \eqref{eq:tildenablaK-transition}
we obtain
\begin{eqnarray}
&&
\tilde{\nabla}_{b} \tilde{K}_{a_{1} a_{2} \ldots a_{p}} +
\sum_{j=1}^p M^{c_{j}}{}_{b a_{j}} (\tilde{\boldsymbol\Lambda},\tilde{\boldsymbol g},p,s) \Theta^{s}
K_{a_{1} \ldots c_{j} \ldots a_{p}} =\nonumber
\\
&&
\Theta^{s} \left[ \nabla_{b} K_{a_{1} \ldots a_{p}} +
\sum_{j=1}^p M^{c_{j}}{}_{b a_{j}}
({\boldsymbol\Lambda},{\boldsymbol g},p,s)
K_{a_{1} \ldots c_{j} \ldots a_{p}} \right]
\label{eq:conformal-diff}
\end{eqnarray}

The previous relation can be cast
in terms of a new conformally covariant operator
$D_{b}^{s,(0,p)}:\mathcal{T}^0_p(\mathcal{M})\rightarrow\mathcal{T}^0_{p+1}(\mathcal{M})$:

\begin{equation}
\tilde{D}_{b}^{s,(0,p)} \tilde{K}_{a_{1} \ldots a_{p}} =
\Theta^{s} D_{b}^{s, (0, p)} K_{a_{1} \ldots a_{p}}.
\label{eq:conformal-Dsp}
\end{equation}

Here, the operator ${D}_{b}^{s,(0, p)}$ is defined by (we
shorten the notation by suppressing the covariant indices of $K$
which are those of \eqref{eq:conformal-diff})
\begin{eqnarray}
&&
D_{b}^{s, (0,p)}\cdot K \equiv \nabla_{b}\cdot K
+ \sum_{j=1}^p M^{c}{}_{b a}({\boldsymbol\Lambda},{\boldsymbol g},p,s)\cdot K =\\
&&
\nabla_{b} \cdot K + \sum_{j=1}^p\left(\left( \frac{s + p}{p} \right) \Lambda_{b} \delta_{a}{}^{c} +
\Lambda_{a} \delta_{b}{}^{c} - g_{a b} g^{cd}\Lambda_{d}\right)\cdot K.
\label{eq:define-Dsp}
\end{eqnarray}

Equation \eqref{eq:conformal-Dsp} shows that the operator ${D}_{b}^{s, p}$
is conformally covariant for conformally covariant tensors of rank $p$
and conformal weight $s$.

In a similar fashion we can carry out the computations
for conformally contravariant tensorial
metric concomitants. The result is
\[
\tilde{\nabla}_{b} \tilde{K}^{a_{1} a_{2} \ldots a_{p}} +
\sum_{j} \bar M^{a_{j}}{}_{b c_{j}} (\tilde{\boldsymbol\Lambda},\tilde{\boldsymbol g},p,s) \Theta^{s}
K^{a_{1} \ldots c_{j} \ldots a_{p}} =
\]

\[
\Theta^{s} \left[ \nabla_{b} K^{a_{1} \ldots a_{p}} +
\sum_{j=1}^p\bar M^{a_{j}}{}_{b c_{j}} ({\boldsymbol\Lambda},{\boldsymbol g})
K^{a_{1} \ldots c_{j} \ldots a_{p}} \right],
\]
where
\begin{equation}
 \bar M^c{}_{ba}({\boldsymbol\Lambda},{\boldsymbol g},p,s)
 \equiv \left( \frac{s - p}{p} \right) \Lambda_{b} \delta_{a}{}^{c} - \Lambda_{a} \delta_{b}{}^{c}  + g_{a b} g^{c r} \Lambda_{r}
\end{equation}

Again, the previous relation can be cast
in terms of a new conformally covariant operator
$D_{b}^{s, (p, 0)}:\mathcal{T}_0^p(\mathcal{M})\rightarrow\mathcal{T}_0^{p+1}(\mathcal{M})$:

\begin{equation}
\tilde{D}_{b}^{s, (p, 0)} \tilde{K}^{a_{1} \ldots a_{p}} =
\Theta^{s} D_{b}^{s, (p, 0)} K^{a_{1} \ldots a_{p}}.
\label{eq:conformal-barDsp}
\end{equation}

Here, the operator $\bar{D}_{b}^{s, (p,0)}$ is defined by
(again we suppress irrelevant contravariant indices to shorten the notation)
\begin{eqnarray}
&&
D_{b}^{s, (p,0)}\cdot K \equiv \nabla_{b}\cdot K
+ \sum_{j=1}^p\bar M^{c}{}_{b a}({\boldsymbol\Lambda},{\boldsymbol g},p,s)\cdot K =\nonumber\\
&&
\nabla_{b}\cdot K + \sum_{j=1}^p \left(\left( \frac{s - p}{p} \right) \Lambda_{b} \delta_{a}{}^{c} -
\Lambda_{a} \delta_{b}{}^{c} + g_{a b} g^{cd}\Lambda_{d}\right)\cdot K.
\label{eq:define-barDsp}
\end{eqnarray}

Equation \eqref{eq:conformal-barDsp} shows that the operator $\bar{D}_{b}^{s, (p,0)}$
is conformally covariant for conformally contravariant tensors of rank $p$
and conformal weight $s$.

Finally, we can also extend these notions to mixed conformally covariant tensors.
For example, if $K^{p}{}_{q}$ is a metric concomitant conformally covariant
\begin{equation}
 \tilde{K}^{p}{}_{q}=\Theta^s K^p{}_q,
\end{equation}
then
\begin{equation}
\tilde{D}_b^{s,(1,1)} \tilde{K}^{p}{}_{q} = \Theta^s
D_b^{s,1,1} K^{p}{}_{q},
\end{equation}
where
\begin{equation}
 D_b^{s,(1,1)} K^{p}{}_{q} = \nabla_b K^{p}{}_{q} + \bar{M}^p{}_{br}({\boldsymbol\Lambda},{\boldsymbol g}) K^{r}{}_{q} +
 M^r{}_{bq}({\boldsymbol\Lambda},{\boldsymbol g}) K^p{}_r.
\end{equation}
We summarize our previous findings in the following Theorem.

\begin{theorem}[Conformal pseudo-differential]
 Let the tensor $K^{a_1\dots a_p}{}_{b_1\dots b_q}$ be an unphysical metric
 concomitant and assume that it is conformally covariant
 of weight $s\in\mathbb{Z}$
 \begin{equation}
  \widetilde K^{a_1\dots a_p}{}_{b_1\dots b_q}=
  \Theta^s K^{a_1\dots a_p}{}_{b_1\dots b_q}.
 \end{equation}
Then its \underline{conformal pseudo-differential}
$D_a^{s,(p,q)}:\mathcal{T}^p_q(\mathcal{M})\rightarrow \mathcal{T}^p_{q+1}(\mathcal{M})$,
$(p,q)\neq (0,0)$ defined by
\begin{eqnarray}
&&
 D_a^{s,(p,q)} K^{a_1\dots a_p}{}_{b_1\dots b_q} \equiv
 \nabla_a K^{a_1\dots a_p}{}_{b_1\dots b_q}+\nonumber\\
&&
\sum_{j=1}^p\bar M^{a_j}{}_{ac}({\boldsymbol\Lambda},{\boldsymbol g},p,s) K^{a_1\dots a_{j-1}c a_{j+1}\dots a_p}{}_{b_1\dots b_q}+\nonumber\\
&&
\sum_{j=1}^q M^{c}{}_{ab_j}({\boldsymbol\Lambda},{\boldsymbol g},p,s)K^{a_1\dots a_p}{}_{b_1\dots b_{j-1}cb_{j+1}\dots b_q}.
\label{eq:dspq-diff}
\end{eqnarray}
is also conformally covariant with the same weight
\begin{equation}
 \widetilde D_a^{s,(p,q)}\tilde K^{a_1\dots a_p}{}_{b_1\dots b_q} =
 \Theta^s  D_a^{s,(p,q)} K^{a_1\dots a_p}{}_{b_1\dots b_q},
\label{eq:dspq-diff-conformal-cov}
\end{equation}
where
\begin{eqnarray}
&&
M^{c}{}_{b a}({\boldsymbol\Lambda},{\boldsymbol g},p,s)\equiv
\left( \frac{s + p}{p} \right) \Lambda_{b} \delta_{a}{}^{c} +
\Lambda_{a} \delta_{b}{}^{c} - g_{a b} g^{cd}\Lambda_{d}\label{eq:define-M},\\
&&
  \bar M^c{}_{ba}({\boldsymbol\Lambda},{\boldsymbol g},p,s)
 \equiv \left( \frac{s - p}{p} \right) \Lambda_{b} \delta_{a}{}^{c} - \Lambda_{a} \delta_{b}{}^{c}  + g_{a b} g^{c r} \Lambda_{r}.\label{eq:define-barM}
\end{eqnarray}
\label{theo::dspq-diff}
\end{theorem}
\qed

If $(p,q)=(0,0)$ (scalar quantity) then we make the definition
\begin{equation}
 D_a^{s,(0,0)}\equiv D^s_a,
\end{equation}
where $D^s_a$ is the operator defined in \eqref{eq:define-conformal-diff-scalars}.

\subsection{The $\mathcal{C}$-connection}
\label{subsec:C-connection}
Note that when $s=0$ (conformal invariance)
the operators $D^{0,(0,p)}_a$, $\bar D^{0,(p,0)}_a$ are actually independent from
$p$ as it follows from \eqref{eq:define-M}-\eqref{eq:define-barM}.
Therefore we can make the following definitions
\begin{equation}
\mathcal{C}_a\equiv D_{a}^{0,(0,p)},\quad
\bar{\mathcal{C}}_a\equiv D_{a}^{0,(p,0)}.
\end{equation}
These operators are true covariant derivatives as one can easily
check from \eqref{eq:dspq-diff}. In fact, we also deduce from \eqref{eq:dspq-diff}
that $\bar{\mathcal{C}}_a=\mathcal{C}_a$, so we only need to describe one
of them (we pick $\mathcal{C}_a$). Using eq. \eqref{eq:define-Dsp} we have that
\begin{equation}
\Gamma[\mathcal{C},\nabla]^c{}_{ab} =
g^{cd}g_{ab}\Lambda_d -\delta_b{}^c \Lambda_a-\delta_a{}^c \Lambda_b.
\label{eq:transition-C-connection}
\end{equation}
We shall denote the Riemann and Ricci tensors defined by the $C$-connection
as $\mathcal{R}_{abc}{}^d$, $\mathcal{R}_{ab}\equiv\mathcal{R}_{adb}{}^d$
respectively. Given that $\nabla$ has no torsion and $\Gamma[\mathcal{C},\nabla]^c{}_{ab}$
is symmetric, the connection $\mathcal{C}$ has no torsion either.
Hence, the standard formula that relates the Riemann curvature tensors of
$\nabla$  and $\mathcal{C}$ is
\begin{equation}
\mathcal{R}_{abc}{}^{d} =
2 \Gamma [\mathcal{C}, \nabla]^{d}{}_{e[a} \Gamma [\mathcal{C},\nabla]^{e}{}_{b]c}
+ R_{abc}{}^{d} -
2\mathcal{C}_{[a}\Gamma [\mathcal{C},\nabla]^{d}{}_{b]c}.
\label{eq:riemann-formula}
\end{equation}
Inserting \eqref{eq:transition-C-connection} in \eqref{eq:riemann-formula} we get
\begin{eqnarray}
&&
\mathcal{R}_{abc}{}^{d} = R_{abc}{}^{d}-2 \Lambda^{d} \Lambda_{[b}  g_{a]c}+
  2\delta_{[a}{}^{d} (\Lambda_{b]} \Lambda_{c} -  g_{b]c} \Lambda_{e} \Lambda^{e} -  \mathcal{C}_{b]}\Lambda_{c})+\nonumber\\
&&
  2\delta_{c}{}^{d}\mathcal{C}_{[a}\Lambda_{b]}  - 2 g^{de} g_{c[b} \mathcal{C}_{a]}\Lambda_{e}.
\label{eq:riemann-to-mathcalriemann}
\end{eqnarray}

In terms of the $\mathcal{C}$-connection, equation \eqref{eq:dspq-diff-conformal-cov}
takes the form
\begin{equation}
 \widetilde{\mathcal{C}}_a\tilde K^{a_1\dots a_p}{}_{b_1\dots b_q} =
 \mathcal{C}_a K^{a_1\dots a_p}{}_{b_1\dots b_q}.
\label{eq:dspq-diff-mathcal-c}
\end{equation}

Indeed, this result is a consequence of the following stronger one.

\begin{theorem}
 The $\mathcal C$-connection is conformally invariant. This means, using the above notation
\begin{equation}
 \tilde{\mathcal{C}}_a=\mathcal{C}_a.
\end{equation}
\label{theo:c-connection-conf-invariant}
\end{theorem}

\proof
To prove this result we are going to compute the transition tensor
$\Gamma [\mathcal{C} ,\tilde{\mathcal{C}}]^{a}{}_{bc}$
and show that it vanishes. To that end we start from the relation
\begin{equation}
\Gamma [\mathcal{C} ,\tilde{\mathcal{C}}]^{a}{}_{bc} =
\Gamma[\mathcal{C} ,\tilde{\nabla}]^{a}{}_{bc}
+  \Gamma [\tilde{\nabla},\tilde{\mathcal{C}}]^{a}{}_{bc}.
\end{equation}
The last term of the right hand side can be worked out by writing \eqref{eq:transition-C-connection}
in terms of the physical metric, whereas the first term comes from the relation
\begin{equation}
 \Gamma[\mathcal{C} ,\tilde{\nabla}]^{a}{}_{bc} =
\Gamma [\mathcal{C},\nabla]^{a}{}_{bc}+\Gamma [\nabla,\tilde{\nabla}]^{a}{}_{bc}.
\end{equation}
The values of the transition tensors $\Gamma [\nabla ,\tilde{\nabla}]^{a}{}_{bc}$
and $\Gamma [\mathcal{C}, \nabla ]^{a}{}_{bc}$ are given by respectively eqns.
\eqref{eq:define-connection-tensor} and \eqref{eq:transition-C-connection}.
By doing all these replacements we get
\begin{equation}
\Gamma [\mathcal{C} ,\tilde{\mathcal{C}}]^{a}{}_{bc} = \
2\delta_{(c}{}^{a} \tilde{\Lambda}_{b)}
-  \tilde{g}^{ad}\tilde{g}_{bc}\tilde{\Lambda}_{d}
-2\delta_{(c}{}^{a} \Lambda_{b)} + \Lambda^{a} g_{cb}
+ 2 \delta_{(c}{}^{a} \Upsilon_{b)} - g^{ad} g_{bc} \Upsilon_{d}.
\end{equation}
The result follows after using \eqref{eq:unphysicaltophysical-downstairs},
\eqref{eq:unphysicaltophysical-upstairs} and Lemma \ref{lemm:lambda-upsilon}
in the previous equation.
\qed

\begin{remark}\em
A direct consequence of Theorem \ref{theo:c-connection-conf-invariant}
is the conformal invariance of the curvature of the $\mathcal{C}$-connection
\begin{equation}
 \tilde{\mathcal{R}}_{abc}{}^d=
 \mathcal{R}_{abc}{}^d.
\end{equation}
From this we deduce that $\mathcal{R}_{ab}$ is also conformally invariant.
\end{remark}

\section{Conformal Einstein spaces  in terms of the $\mathcal{C}$-connection}
\label{sec:conformal-einstein}

The following result was originally presented in \cite{GOVER2006450}
and reformulated in \cite{GARCIAPARRADO2024105093}. We choose
the latter form as it is better adapted to our purposes.

\begin{theorem}[Gover \& Nurowski 2006]
 Let $(\mathcal{M},g_{ab})$ be a 4-dimensional spacetime and define the following concomitants of $g_{ab}$.
 Define the unphysical metric concomitants
\begin{equation*}
\Lambda^d\equiv\left(\frac{8}{C\cdot C}\right)C^{damp}\nabla_{[m}L_{p]a},\quad
C\cdot C\equiv C_{abcd}C^{abcd}.
\end{equation*}
Then, assuming that $C\cdot C\neq 0$, $(\mathcal{M},g_{ab})$ is locally conformal to an Einstein spacetime if and only if the following
covariant concomitant of $g_{ab}$ vanishes
\begin{equation}
 E_{ab}(g_{ab},\nabla)\equiv\mbox{\rm TF}_{g}\big[
	L_{ab}+
	\nabla_{a}\Lambda_{b}
	+ \Lambda_a\Lambda_b\big],
\label{eq:ideal-einstein}
\end{equation}
where ${\rm TF}_g$ stands for ``trace-free part'' with respect to $g_{ab}$.
Moreover $E_{ab}$ is \emph{conformally invariant} according to
\begin{equation*}
 \tilde{E}_{ab}(\tilde{g}_{ab},\tilde\nabla)= E_{ab}(g_{ab},\nabla).
\end{equation*}
\label{theo:ideal-einstein}
\end{theorem}
\proof See \cite{GARCIAPARRADO2024105093}.\qed

Using the previous theorem as the starting point
we obtain the main result of the present work.
\begin{theorem}
 A 4-dimensional
 manifold $(\mathcal{M}, g)$ is conformally related to an Einstein manifold
 $(\tilde{\mathcal{M}},\tilde{g}_{ab})$, if and only if its $\mathcal{C}$-Ricci
 tensor fulfills the following algebraic conditions
 \begin{equation}
\mathcal{R}_{[ab]}=0,\quad
\mathcal{R}_{ab}-\frac{1}{4} g_{ab} \mathcal{R} =0.
\label{eq:C-Ricci-Einstein}
\end{equation}
\end{theorem}

\proof
Using \eqref{eq:riemann-to-mathcalriemann}
one can find the following relations
\begin{eqnarray}
&&
R_{ac} = 2 \Lambda_{a} \Lambda_{c} + \mathcal{R}_{ac} -
2(\Lambda_{b}\Lambda^{b}) g_{ac} - 3 \mathcal{C}_{a}\Lambda_{c} + \
\mathcal{C}_{c}\Lambda_{a} -  g_{ac} (g^{db} \mathcal{C}_{d}\Lambda_{b}),
\label{eq:c-ricci}\\
&&
R = -6 \Lambda_{a} \Lambda^{a} + \mathcal{R}_{a}{}^{a} - 6 g^{ab} \
\mathcal{C}_{b}\Lambda_{a}
\end{eqnarray}

From this we can obtain the relation between the Schouten tensor
$L_{ab}$ and $\mathcal{R}_{ab}$ by means of \eqref{eq:define-schouten}
\begin{equation}
L_{ab} = \Lambda_{a} \Lambda_{b} -  \frac{1}{2} \Lambda_{c} \Lambda^{c} \
g_{ab} +
\frac{1}{2} \mathcal{R}_{(ab)}
-  \frac{1}{12} \mathcal{R}_{c}{}^{c} g_{ab} - \mathcal{C}_{(a}\Lambda_{b)}.
\end{equation}
Next we use the relation between the $\nabla$ covariant derivative
and the $\mathcal{C}$ covariant derivative
\begin{equation}
\nabla_{b}\Lambda_{a} = - \Gamma [\nabla ,\mathcal{C} ]^{c}{}_{ba} \
\Lambda_{c} + \mathcal{C}_{b}\Lambda_{a}.
\end{equation}
Since the connection tensor $\Gamma [\nabla ,\mathcal{C} ]^{c}{}_{ba}$
is known (see \eqref{eq:transition-C-connection}) we can replace it
getting
\begin{equation}
 \nabla_{b}\Lambda_{a} =
 -2 \Lambda_{a} \Lambda_{b} + \Lambda_{c} \Lambda^{c} g_{ab} + \
\mathcal{C}_{b}\Lambda_{a}.
\end{equation}
Putting together all the pieces we can compute $E_{ab}$ as defined by
\eqref{eq:ideal-einstein} but written now in terms of the
$\mathcal{C}$-connection. The result is
\begin{equation}
 E_{ab}(g_{ab},\nabla)=\mbox{\rm TF}_{g}\left[
\frac{1}{2} \mathcal{R}_{(ab)}
+ (\tfrac{1}{2} \Lambda_{c} \Lambda^{c} -  \tfrac{1}{12} \
\mathcal{R}_{c}{}^{c}) g_{ab} - \mathcal{C}_{[a}\Lambda_{b]}
\right].
\end{equation}
From \eqref{eq:c-ricci} we deduce
\begin{equation}
 \mathcal{C}_{[a}\Lambda_{b]} = \tfrac{1}{4} \mathcal{R}_{[ab]}.
\end{equation}
Thus
\begin{equation}
 E_{ab}(g_{ab},\nabla)=\mbox{\rm TF}_{g}\left[
\frac{1}{2} \mathcal{R}_{(ab)}\right]
 - \frac{1}{4} \mathcal{R}_{[ab]}.
\end{equation}
Hence $E_{ab}$ vanishes if and only if
\begin{equation}
\mathcal{R}_{[ab]}=0,\quad
\mathcal{R}_{ab}-\frac{1}{4} g_{ab} \mathcal{R} =0,
\end{equation}
as desired.

\qed

\begin{remark}\em
 In the case that $(\mathcal{M},g_{ab})$ is a \emph{Lorentzian manifold}
 (space-time) we shall refer to the equation set \eqref{eq:C-Ricci-Einstein}
 as the \emph{conformal vacuum Einstein equations}.

\end{remark}


\section{Conclusions}
In this work, we have studied the structure of a 4-dimensional pseudo-Riemannian manifold in the context where
$C_{abcd}C^{abcd} \neq 0$. A key outcome of our research is the introduction of the \emph{conformal pseudo-
differential}, a tool that enables the construction of new tensorial conformally covariant metric concomitants
from pre-existing ones.

Building on the conformal pseudo-differential, we defined the $\mathcal{C}-$connection, a specialization of the
Weyl connection. This connection possesses the notable property that its Riemann tensor remains conformally
invariant, offering a new approach to understanding the conformal structure of such manifolds. Furthermore, we
presented a novel characterization of Einstein spaces through the $\mathcal{C}-$connection, as detailed in Theorem
\ref{theo:ideal-einstein}, providing an alternative perspective on their geometric structure.

A particularly intriguing open problem emerging from this work is whether established conformal operators, such as
the Yamabe or Paneitz operators, can be reformulated using the conformal pseudo-differential that we have
introduced. Another interesting question is the formulation of an initial value
problem for the conformal vacuum Einstein equations \eqref{eq:C-Ricci-Einstein}.
These problems will be addressed elsewhere.

\section*{Acknowledgments}
\noindent
Supported by Spanish MICINN Project No. PID2021-126217NB-I00.

\providecommand{\bysame}{\leavevmode\hbox to3em{\hrulefill}\thinspace}
\providecommand{\MR}{\relax\ifhmode\unskip\space\fi MR }
\providecommand{\MRhref}[2]{%
  \href{http://www.ams.org/mathscinet-getitem?mr=#1}{#2}
}
\providecommand{\href}[2]{#2}

\end{document}